# A homotopy method for finding eigenvalues and eigenvectors

Kerry M. Soileau

August 29, 2007

ABSTRACT

Suppose we want to find the eigenvalues and eigenvectors for the linear operator $L$, and suppose that we have solved this problem for some other "nearby" operator $K$. In this paper we show how to represent the eigenvalues and eigenvectors of $L$ in terms of the corresponding properties of $K$.

## INTRODUCTION

Suppose we want to find the eigenvalues and eigenvectors for the linear operator $L$ over some Hilbert space with inner product $\langle \cdot, \cdot \rangle$. That is, we want to find the vectors $x$ and the scalars $\lambda$ satisfying $Lx = \lambda x$. Also suppose that we have solved this problem for some other linear operator $K$. Is it possible to represent the eigenvalues and eigenvectors of $L$ in terms of the corresponding properties of $K$? In this paper we address this question.

## APPROACH

We use a homotopy technique. We form a convex combination $\theta L + (1-\theta) K$ of the operators $K$ and $L$ with parameter $\theta$. As the parameter ranges from $0$ to $1$, the operator $\theta L + (1-\theta) K$ ranges from $K$ to $L$. More precisely, at each value of $\theta$ the corresponding operator is given by $(\theta L + (1-\theta) K) x = \theta L x + (1-\theta) K x$

Next, we envision, for each choice of $\theta$, a set of eigenvectors $\{x_i(\theta)\}_{i=1}^{\infty}$ and eigenvalues $\{\lambda_i(\theta)\}_{i=1}^{\infty}$ for the corresponding operator $\theta L + (1-\theta) K$, satisfying

$\lambda_i(\theta) x_i(\theta) = \theta L x_i(\theta) + (1-\theta) K x_i(\theta)$. Then $x_i(0) = e_i$, where $\{e_i\}_{i=1}^{\infty}$ are eigenvectors of $K$, and $\{\lambda_i(0)\}_{i=1}^{\infty}$ are the corresponding eigenvalues.

Our strategy is to develop power series representations of $\{x_i(\theta)\}_{i=1}^{\infty}$ and $\{\lambda_i(\theta)\}_{i=1}^{\infty}$, and then, if possible, to evaluate them at the value $\theta = 1$. If the power series representations converge absolutely there, we may expect to find that $\{x_i(1)\}_{i=1}^{\infty}$ and $\{\lambda_i(1)\}_{i=1}^{\infty}$ are eigenvectors and eigenvalues of the operator $L$.

**DERIVATION**

We stipulated earlier that $\{e_i\}_{i=1}^{\infty}$ are eigenvectors of $K$, whence $x_i(0) = e_i$. We also assume that they are orthonormal: $\langle e_m, e_n \rangle = \delta_m^n$. Let $\{\lambda_i(0)\}_{i=1}^{\infty}$ represent their corresponding eigenvalues, so that $\lambda_i(0) e_i = K e_i = K x_i(0)$. We immediately get $\langle K e_i, e_j \rangle = \delta_i^j \lambda_i(0)$. In case $i = j$ we get

$$\boxed{\lambda_i(0) = \langle K e_i, e_i \rangle}$$

Next, we define $\{\varphi_{i,j}(\theta)\}_{i,j=1}^{\infty}$ as $\varphi_{i,j}(\theta) = \langle x_i(\theta), e_j \rangle$. We then have that

$x_i(\theta) = \sum_{j=1}^{\infty} \varphi_{i,j}(\theta) e_j$. In particular $\varphi_{i,j}(0) = \langle x_i(0), e_j \rangle = \langle e_i, e_j \rangle = \delta_i^j$, whence

$$\boxed{\varphi_{i,j}(0) = \delta_j^i}$$

Making the appropriate substitutions into $\lambda_i(\theta) x_i(\theta) = (\theta L + (1-\theta) K) x_i(\theta)$ yields

$$\lambda_i(\theta) \sum_{j=1}^{\infty} \varphi_{i,j}(\theta) e_j = (\theta L + (1-\theta) K) \sum_{j=1}^{\infty} \varphi_{i,j}(\theta) e_j.$$

After rearrangement we get $\lambda_i(\theta)\sum_{j=1}^{\infty}\varphi_{i,j}(\theta)e_j = \theta\sum_{j=1}^{\infty}\varphi_{i,j}(\theta)Le_j + (1-\theta)\sum_{j=1}^{\infty}\varphi_{i,j}(\theta)Ke_j$

from which there results the equation

$$\lambda_i(\theta)\sum_{j=1}^{\infty}\varphi_{i,j}(\theta)\langle e_j, e_k\rangle = \theta\sum_{j=1}^{\infty}\varphi_{i,j}(\theta)\langle Le_j, e_k\rangle + (1-\theta)\sum_{j=1}^{\infty}\varphi_{i,j}(\theta)\langle Ke_j, e_k\rangle$$

Recalling $\langle e_i, e_j\rangle = \delta_i^j$ and defining

$X_{m,n} \equiv \langle Le_m, e_n\rangle$ leads, after substitution, to the equation

$$\lambda_i(\theta)\sum_{j=1}^{\infty}\varphi_{i,j}(\theta)\delta_j^k = \theta\sum_{j=1}^{\infty}\varphi_{i,j}(\theta)X_{j,k} + (1-\theta)\sum_{j=1}^{\infty}\varphi_{i,j}(\theta)\delta_j^k\lambda_j(0)$$

Simplifying yields

$$\boxed{\lambda_i(\theta)\varphi_{i,k}(\theta) = \theta\sum_{j=1}^{\infty}\varphi_{i,j}(\theta)X_{j,k} + (1-\theta)\varphi_{i,k}(\theta)\lambda_k(0)}$$

The approach from this point is to calculate about $\theta = 0$ the Maclaurin series expansions of the functions $\{x_i(\theta)\}_{i=1}^{\infty}$ and $\{\lambda_i(\theta)\}_{i=1}^{\infty}$. To do this our procedure is to solve for the coefficients of power series which formally solve the above equation. That is, we assume the forms $\lambda_i(\theta) = \sum_{n=0}^{\infty} a_{i,n}\theta^n$ and $\varphi_{i,j}(\theta) = \sum_{r=0}^{\infty} b_{i,j,r}\theta^r$ and substitute them in that equation, then solve for the coefficients $\{a_{i,n}\}$ and $\{b_{i,j,r}\}$.

After performing the substitution, we get

$\sum_{n=0}^{\infty} a_{i,n}\theta^n \sum_{m=0}^{\infty} b_{i,k,m}\theta^m = \theta\sum_{j=1}^{\infty}\sum_{r=0}^{\infty} b_{i,j,r}\theta^r X_{j,k} + (1-\theta)\sum_{r=0}^{\infty} b_{i,k,r}\theta^r \lambda_k(0)$. Upon rearranging, there results $\sum_{s=0}^{\infty}\theta^s \sum_{m=0}^{s} a_{i,s-m}b_{i,k,m} = \sum_{r=1}^{\infty}\theta^r \sum_{j=1}^{\infty} X_{j,k}b_{i,j,r-1} + \sum_{r=0}^{\infty}\lambda_k(0)b_{i,k,r}\theta^r + \sum_{r=1}^{\infty} -\lambda_k(0)b_{i,k,r-1}\theta^r$. We first set equal the coefficients of the constant terms to get $\sum_{m=0}^{0} a_{i,0-m} b_{i,k,m} = \lambda_k(0) b_{i,k,0}$, which when simplified becomes $\left(a_{i,0} - \lambda_k(0)\right) b_{i,k,0} = 0$. Since $\lambda_i(\theta) = \sum_{n=0}^{\infty} a_{i,n} \theta^n$, we must have the equation $a_{i,0} = \lambda_i(0)$. Since $\varphi_{i,j}(\theta) = \sum_{r=0}^{\infty} b_{i,j,r} \theta^r$, we must have the equation

$b_{i,j,0} = \varphi_{i,j}(0) = \delta_j^i$, and these two equations indeed satisfy $\left(a_{i,0} - \lambda_k(0)\right) b_{i,k,0} = 0$. Next,

for $r > 0$ we set equal the coefficients of the $\theta^r$ terms to get

$\sum_{m=0}^{r} a_{i,r-m} b_{i,k,m} = \sum_{j=1}^{\infty} X_{j,k} b_{i,j,r-1} + \lambda_k(0) b_{i,k,r} - \lambda_k(0) b_{i,k,r-1}$. Upon substituting $b_{i,k,0} = \delta_i^k$ we

get $a_{i,r} \delta_i^k + \sum_{m=1}^{r} a_{i,r-m} b_{i,k,m} = \sum_{j=1}^{\infty} X_{j,k} b_{i,j,r-1} + \lambda_k(0) b_{i,k,r} - \lambda_k(0) b_{i,k,r-1}$. The examination of

this equation splits naturally into the two cases $i = k$ and $i \neq k$. Considering the case

$i = k$ first, we get $a_{i,r} = \sum_{j=1}^{\infty} X_{j,i} b_{i,j,r-1} + \lambda_i(0) b_{i,i,r} - \lambda_i(0) b_{i,i,r-1} - \sum_{m=1}^{r} a_{i,r-m} b_{i,i,m}$. In case

$i \neq k$, we get $\sum_{m=1}^{r} a_{i,r-m} b_{i,k,m} = \sum_{j=1}^{\infty} X_{j,k} b_{i,j,r-1} + \lambda_k(0) b_{i,k,r} - \lambda_k(0) b_{i,k,r-1}$, which may be

rearranged as $a_{i,0} b_{i,k,r} = \sum_{j=1}^{\infty} X_{j,k} b_{i,j,r-1} + \lambda_k(0) \left(b_{i,k,r} - b_{i,k,r-1}\right) - \sum_{m=1}^{r-1} a_{i,r-m} b_{i,k,m}$. After

substituting $a_{i,0} = \lambda_i(0)$ and solving, we get

$$b_{i,k,r} = \frac{1 - \delta_i^k}{\lambda_i(0) - \lambda_k(0)} \left( \sum_{j=1}^{\infty} X_{j,k} b_{i,j,r-1} - \sum_{m=1}^{r-1} a_{i,r-m} b_{i,k,m} - \lambda_k(0) b_{i,k,r-1} \right)$$

We now perform this recursion to get some initial results. Using

$a_{i,r} = \sum_{j=1}^{\infty} X_{j,i} b_{i,j,r-1} + \lambda_i(0) b_{i,i,r} - \lambda_i(0) b_{i,i,r-1} - \sum_{m=1}^{r} a_{i,r-m} b_{i,i,m}$, we see that

$$a_{i,1} = \sum_{j=1}^{\infty} X_{j,i} b_{i,j,0} + \lambda_i(0) b_{i,i,1} - \lambda_i(0) b_{i,i,0} - \sum_{m=1}^{1} a_{i,1-m} b_{i,i,m} \text{ hence}$$

$$a_{i,1} = \sum_{j=1}^{\infty} X_{j,i} \delta_i^j + \lambda_i(0) b_{i,i,1} - \lambda_i(0) \delta_i^i - \lambda_i(0) b_{i,i,1} \text{ and thus } a_{i,1} = X_{i,i} - \lambda_i(0). \text{ Using}$$

$$b_{i,k,r} = \frac{1-\delta_i^k}{\lambda_i(0) - \lambda_k(0)} \left( \sum_{j=1}^{\infty} X_{j,k} b_{i,j,r-1} - \sum_{m=1}^{r-1} a_{i,r-m} b_{i,k,m} - \lambda_k(0) b_{i,k,r-1} \right), \text{ we get that}$$

$$b_{i,k,1} = \frac{1-\delta_i^k}{\lambda_i(0) - \lambda_k(0)} \left( \sum_{j=1}^{\infty} X_{j,k} b_{i,j,0} - \sum_{m=1}^{0} a_{i,r-m} b_{i,k,m} - \lambda_k(0) b_{i,k,0} \right) \text{ and thus}$$

$$b_{i,k,1} = (1-\delta_i^k) \frac{X_{i,k}}{\lambda_i(0) - \lambda_k(0)}. \text{ Similarly we use the recursion formula to produce}$$

$$a_{i,2} = \sum_{j=1}^{\infty} X_{j,i} b_{i,j,1} + \lambda_i(0)(b_{i,i,2} - b_{i,i,1}) - \sum_{m=1}^{2} a_{i,2-m} b_{i,i,m}, \text{ which after substitution yields}$$

$$a_{i,2} = \sum_{j=1}^{\infty} X_{j,i} (1-\delta_i^j) \frac{X_{i,j}}{\lambda_i(0) - \lambda_j(0)} = \sum_{\substack{j=1 \\ j \neq i}}^{\infty} \frac{X_{i,j} X_{j,i}}{\lambda_i(0) - \lambda_j(0)}. \text{ Finally, we get}$$

$$b_{i,k,2} = \frac{1-\delta_i^k}{\lambda_i(0) - \lambda_k(0)} \left( \sum_{j=1}^{\infty} X_{j,k} b_{i,j,1} - \sum_{m=1}^{1} a_{i,2-m} b_{i,k,m} - \lambda_k(0) b_{i,k,1} \right) \text{ which after substitution}$$

becomes $b_{i,k,2} = \frac{1-\delta_i^k}{\lambda_i(0) - \lambda_k(0)} \left( \sum_{\substack{j=1 \\ j \neq i}}^{\infty} X_{j,k} \frac{X_{i,j}}{\lambda_i(0) - \lambda_j(0)} + X_{i,k} \left( 1 - \frac{X_{i,i}}{\lambda_i(0) - \lambda_k(0)} \right) \right)$. This

means that

$$\lambda_i(\theta) = \lambda_i(0) + \left(X_{i,i} - \lambda_i(0)\right)\theta + \sum_{\substack{j=1 \\ j \neq i}}^{\infty} \frac{X_{i,j} X_{j,i}}{\lambda_i(0) - \lambda_j(0)} \theta^2 + \cdots$$ and

$$x_i(\theta) = e_i + \sum_{\substack{j=1 \\ j \neq i}}^{\infty} \left( \frac{X_{i,j}}{\lambda_i(0) - \lambda_j(0)} \theta + \left( \frac{1}{\lambda_i(0) - \lambda_j(0)} \sum_{\substack{k=1 \\ k \neq i}}^{\infty} \frac{X_{i,k} X_{k,j}}{\lambda_i(0) - \lambda_k(0)} + \frac{X_{i,j}}{\lambda_i(0) - \lambda_j(0)} - \frac{X_{i,i} X_{i,j}}{\left(\lambda_i(0) - \lambda_j(0)\right)^2} \right) \theta^2 + \cdots \right) e_j$$

(at least formally). Now if it happens that these formal power series actually converge absolutely at $\theta = 1$, we then have $$\lambda_i(1) = X_{i,i} + \sum_{\substack{j=1 \\ j \neq i}}^{\infty} \frac{X_{i,j} X_{j,i}}{\lambda_i(0) - \lambda_j(0)} + \cdots$$ and

$$x_i(1) = e_i + \sum_{\substack{j=1 \\ j \neq i}}^{\infty} \left( \frac{X_{i,j}}{\lambda_i(0) - \lambda_j(0)} + \left( \frac{1}{\lambda_i(0) - \lambda_j(0)} \sum_{\substack{n=1 \\ n \neq i}}^{\infty} \frac{X_{i,n} X_{n,j}}{\lambda_i(0) - \lambda_n(0)} + \frac{X_{i,j}}{\lambda_i(0) - \lambda_j(0)} - \frac{X_{i,i} X_{i,j}}{\left(\lambda_i(0) - \lambda_j(0)\right)^2} \right) + \cdots \right) e_j.$$

These are eigenvalues and eigenvectors of the operator $L$, as desired.